\documentclass[12pt]{amsart}

\usepackage{amssymb,amsmath}
\usepackage{url}
\usepackage{amscd}
\usepackage{texdraw}
\usepackage[all]{xy}
\usepackage{fancyhdr}
\usepackage{amsmath}
\usepackage{graphicx}
\usepackage{subcaption}

\usepackage{amscd,color}
\usepackage[shortlabels]{enumitem}

\usepackage{amsfonts}
\usepackage{amssymb}
\usepackage{amsmath}
\usepackage{amsthm}

\usepackage{float}

\normalsize

\newtheorem{theorem}{Theorem}
\newtheorem{example}{Example}

\newtheorem{lemma}{Lemma}
\newtheorem{proposition}{Proposition}
\newtheorem{problem}{Problem}
\newtheorem{corollary}{Corollary}
\newtheorem{conjecture}{Conjecture}

\theoremstyle{definition}

\newcommand{\C}{\mathbb C}

\newcommand{\T}{\mathbb{T}}

\begin{document}

\title[Geometric Mean Value]{Smale's Mean Values And Geometric Mean Values}


\author{Yunping Jiang}
\address[YPJ]{Department of Mathematics\\
 Queens College of the City University of New York\\
 Flushing, NY 11367-1597\\
 and\\
 The Ph.D. Program in Mathematics\\
 Graduate Center of the City University of New York;
 New York, NY 10016}
 \email[Jiang]{yunping.jiang@qc.cuny.edu}

\author{Zhe Wang}
\address[ZW]{Department of Mathematics\\ 
Bronx Community College of the City University of
New York, 2155 University Ave, Bronx, NY 10453}
 \email[Wang]{zhe.wang@bcc.cuny.edu}

 

\thanks{The work is partially supported by PSC-CUNY awards, a PSC-CUNY enhanced award (66680-54), and Simons Foundation awards (523341 and 942077).}


\begin{abstract}
We prove that, for every integer $d>2$ and every $L>0$, there exists a polynomial of degree $d$ whose geometric mean value exceeds $L$. This result demonstrates that the geometric mean value alone is insufficient to establish the smallest upper bound for the minimum of Smale’s mean values. Under certain additional assumptions, we derive upper bounds for the geometric mean value, thereby verifying Smale’s mean value conjecture for several subclasses of polynomials. Furthermore, we connect Smale’s mean values to Fekete points. Motivated by this connection, we formulate a new problem concerning Fekete points in the spirit of Smale’s mean values.
\end{abstract}

\maketitle

\section{\bf \large Introduction}

Let $\C$ be the complex plane. Let ${\mathcal P}={\mathcal P} (d)$ be the space of all complex polynomials $P$ of degree $d\geq 2$, i.e., 
\[
P(z) = a_{d} z^{d} + \cdots + a_{1} z + a_{0}, \quad a_{d} \neq 0, \cdots, a_{1}, a_{0}\in \C. 
\]

Every $P\in {\mathcal P}$ has $d-1$ critical points (counted by multiplicity), which are the roots of its derivative $P'(z) = 0$.  
Let $c_{1}, \ldots, c_{d-1}$ denote these critical points, and let 
\(
v_{1} = P(c_{1}), \ \ldots, \ v_{d-1} = P(c_{d-1})
\)
denote the corresponding critical values (numbered with multiplicities taken into account). In connection with the analysis of the convergence of Newton’s method for calculating the roots of \(P(z)=0\), Smale~\cite{S} considered the following \(d-1\) mean values 
\[
A_{k,P}(z_{0}) = \frac{P(c_{k}) - P(z_{0})}{c_{k} - z_{0}} \cdot \frac{1}{P'(z_{0})}, 
\quad k = 1, \ldots, d-1,
\]
where $P'(z_{0})\not=0$. By applying the Koebe $1/4$-theorem (refer to, for example,~\cite[Theorem 2.4]{JBook}), we can see directly, 
\[
\min \{|A_{k,P}(z_{0})|\;;\; 1 \leq k \leq d-1\} \leq 4
\]
for all $P\in {\mathcal P}$ (refer to~\cite{SS}). 
It is natural to study the following problem.

\medskip
\begin{problem}~\label{sc}
What is the smallest upper bound $0< M\leq 4$ such that 
\[
\min \{|A_{k,P}(z_{0})|\;;\; 1 \leq k \leq d-1\} \leq M
\]
for all $P\in {\mathcal P}$
\end{problem}

Smale~\cite{SS} conjectured that $M=(d-1)/d$. 
Actually, we know that 
$$
M\leq 4\frac{1+(d-2)4^{\frac{1}{1-d}}}{d+1}
$$ 
(refer to~\cite{CFL,FS}) and for $d\geq 8$, 
$$
M<4-\frac{2.263}{\sqrt{d}}
$$
(refer to ~\cite{C}).  

To study Problem~\ref{sc} and Smale's conjecture, we can consider a normalized polynomial 
\begin{equation}~\label{norm}
Q(z) = \frac{P(az+z_{0})-P(z_{0})}{aP'(z_{0})}, \quad a\not=0, 
\end{equation}
in the sense that $Q(0)=0$ and $Q'(0)=1$. 
The critical points of $Q$ are $\widetilde{c}_{1}=(c_{1}-z_{0})/a$, $\cdots$, $\widetilde{c}_{d-1}=(c_{d-1}-z_{0})/a$ and
$$
A_{k, Q} (0) = A_{k, P} (z_{0}), \quad k=1, \cdots, d-1.
$$
Thus, we need only to consider the subspace 
$$
{\mathcal P}_{1}={\mathcal P}_{1}(d)=\{ P\in {\mathcal P}\;;\;
P(0)=0,\; P'(0)=1\}.
$$ 
For every $P\in {\mathcal P}_{1}$, we still use $c_{1}$, $\cdots$, $c_{d-1}$ to denote its critical points (numbered with multiplicities taken into account). We call the values 
$$
A_{k, P} =A_{k, P}(0)= \frac{P(c_{k})}{c_{k}}, \quad k=1, 2, \cdots, d-1,
$$
Smale's mean values associated with $P\in {\mathcal P}_{1}$. 
Now Problem~\ref{sc} can be restated as 

\medskip
\begin{problem}~\label{scp}
What is the smallest upper bound $M>0$ such that 
\begin{equation}~\label{smv}
\min \{|A_{k,P}|\;;\; 1 \leq k \leq d-1\} \leq M
\end{equation}
for all $P\in {\mathcal P}_{1}$.
\end{problem}

Alongside~(\ref{smv}), we consider the \emph{geometric mean value} associated with $P\in {\mathcal P}_{1}$,
\[
GM_{P} := \left( \prod_{k=1}^{d-1} A_{k,P} \right)^{\frac{1}{d-1}}.
\]
If the inequality
\begin{equation}~\label{meq}
|GM_{P}| \leq M
\end{equation}
were true for all $P\in {\mathcal P}_{1}$, it would immediately imply~(\ref{smv}). It has long been an open question whether we can obtain the best upper bound from the best upper bound for the geometric mean value. 
The main result of this paper answers this question negatively. 

\medskip
\begin{theorem}[Main Theorem]~\label{main}
For any fixed integer $d > 2$ and any real number $L > 0$, there exists a $P\in {\mathcal P}_{1}$ such that
\[
|GM_{P}| > L.
\]
\end{theorem}

Before proving this main result, we first construct a polynomial $P\in {\mathcal P}_{1}$ whose geometric mean value exceeds $(d-1)/d$. We present this in Section~\ref{sec:gmv} as Theorem~\ref{jw}. We prove Theorem~\ref{main} in Section~\ref{sec:pmt}. Theorems~\ref{main} and~\ref{jw} show that the geometric mean value alone cannot be used to study Problem~\ref{scp} and Smale's conjecture. 

However, in Section~\ref{sec:smvc}, we find an upper bound for the geometric mean values for some subclass of polynomials. We present this as Theorem~\ref{unit}. 
Using this upper bound, we reprove Tischler's theorem that Smale's conjecture~\ref{sc} holds for the subclass of polynomials whose roots have the same norm as a corollary (Corollary~\ref{eqroot}). In the same section, we also verify Smale's conjecture for the subclass of polynomials whose critical points are quasi-symmetrically distributed with respect to a circle centered at $0$ (Theorem~\ref{eqcritsymm}). In particular, this implies Smale's conjecture for the subclass of polynomials whose critical points have the same norm (Corollary~\ref{eq3}). 

Also in Section~\ref{sec:smvc}, we show that, by a normalization, we only need to study problem~\ref{scp} for the space of polynomials whose critical points contain $1$. This leads us to connect Problem~\ref{scp} to the study of Fekete points in Section~\ref{sec:smvfp}. This connection lets us reformulate Problem~\ref{scp} as an equivalent problem about the image domain of a several-variable holomorphic function. One advantage of this reformulation is that it provides a framework in which Problem~\ref{scp} can be investigated efficiently with computer assistance. In addition, this equivalent formulation enables us to prove a result concerning Fekete points that is analogous in spirit to Problem~\ref{scp}.
 
\section{\bf \large Geometric Mean Values}~\label{sec:gmv}

For $d=2$ and every $P \in \mathcal{P}_{1}$, 
\[
P(z) = z + a z^{2}, \quad a \in \mathbb{C} \setminus \{0\}.
\]
Since $P'(z) = 1 + 2 a z$, there is only one critical point:
\[
c = -\frac{1}{2a}.
\]
The corresponding critical value is 
\[
v = P(c) = -\frac{1}{2a} + a \cdot \frac{1}{(2a)^{2}} = -\frac{1}{4a}.
\]
Hence
\[
|A_{P}| = |GM_{P}| = \frac{|v|}{|c|} = \frac{1}{2} = \frac{2 - 1}{2}.
\]
Thus, Smale's conjecture holds for all quadratic complex polynomials.

For $d>2$, the following theorem implies that the number $(d-1)/d$ cannot serve as an upper bound for geometric mean values. 

\medskip
\begin{theorem}\label{jw}
For any $d > 2$, there exists a $P \in \mathcal{P}_{1}$ such that 
\[
|GM_{P}| > \frac{d-1}{d}.
\]
\end{theorem}

\begin{proof}
Consider 
\[
P_{0}(z) = z - z^{d} \in \mathcal{P}_{1}.
\]
Since $P_{0}'(z) = 1 - d z^{d-1}$, the critical points of $P_{0}$ are 
\[
c_{k} = \frac{e^{2\pi i k/(d-1)}}{d^{1/(d-1)}}, \quad k = 1, \ldots, d-1,
\]
and the corresponding critical values are
\[
v_{k} = P_{0}(c_{k}) 
= \frac{e^{2\pi i k/(d-1)}}{d^{1/(d-1)}}
- \frac{e^{2\pi i d k/(d-1)}}{d^{d/(d-1)}}
= \frac{e^{2\pi i k/(d-1)}}{d^{1/(d-1)}} \cdot \frac{d-1}{d}
= c_{k} \cdot \frac{d-1}{d}.
\]
Thus, Smale's mean values for $P_{0}$ are
\[
A_{k,P_{0}} = \frac{P_{0}(c_{k})}{c_{k}} = \frac{v_{k}}{c_{k}} = \frac{d-1}{d}, 
\quad k = 1, \ldots, d-1,
\]
and the geometric mean value for $P_{0}$ is  
\[
GM_{P_{0}} = \left( \prod_{k=1}^{d-1} A_{k,P_{0}} \right)^{\frac{1}{d-1}} = \frac{d-1}{d}.
\]

Now let 
\[
\Delta = \{ z \in \mathbb{C} \mid |z| < 1 \}, \qquad
\Delta_{r} = \{ z \in \mathbb{C} \mid |z| < r \}, \ r > 0.
\]
Consider the holomorphic family of polynomials,
\[
\{ P_{\lambda}(z) = z + \lambda z^{d-1} - z^{d}\;|\; \quad \lambda \in \Delta\},
\]
in $\mathcal{P}_{1}$.
Let
\[
c_{1}(\lambda), \ldots, c_{d-1}(\lambda)
\]
be the critical points of $P_{\lambda}$, which depend on $\lambda$ holomorphically. Let 
\[
v_{1}(\lambda) = P_{\lambda}(c_{1}(\lambda)), \ \ldots, \ 
v_{d-1}(\lambda) = P_{\lambda}(c_{d-1}(\lambda))
\]
be the corresponding critical values, which depend on $\lambda$ holomorphically too. For sufficiently small $0 < r < 1$, all $c_{k}(\lambda)$ and $v_{k}(\lambda)$ are nonzero on $\Delta_{r}$.  
We define a holomorphic function
\[
h(\lambda) = GM_{P_{\lambda}} 
= \left( \prod_{k=1}^{d-1} \frac{v_{k}(\lambda)}{c_{k}(\lambda)} \right)^{\frac{1}{d-1}}, 
\quad \lambda \in \Delta_{r}.
\]

Since $h(0) = (d-1)/d$ and since $h$ is an open map, the image $h(\Delta_{r})$ is an open set containing $(d-1)/d$ as an interior point. Therefore, for any $\tau \in h(\Delta_{r})$ with $|\tau| > (d-1)/d$, there exists $\lambda (\tau) \in \Delta_{r}$ such that $h(\lambda(\tau) ) = \tau$.  
Thus for this $P_{\lambda (\tau)} \in \mathcal{P}_{1}$, 
\[
|GM_{P_{\lambda (\tau)}}| > \frac{d-1}{d}.
\]
This completes the proof.
\end{proof}

\section{\bf \large The proof of the main theorem}~\label{sec:pmt}

Before giving our proof of Theorem~\ref{main}. Let us study a general monic complex polynomial 
\[
P(z) = z^{d} + a_{d-1}z^{d-1} + \cdots + a_{2} z^{2} + a_{1} z +a_{0}.
\]
Denote its roots by
\(
z_1, \dots, z_d
\)
and its critical points by
\(
c_1, \dots, c_{d-1}
\)
(numbered with multiplicities taken into account).  
We have the factorizations
\[
P(z) = \prod_{k=1}^d (z - z_k), \qquad 
P'(z) = d \prod_{k=1}^{d-1} (z - c_k).
\]

Differentiating $P(z)$, we have that 
\[
P'(z) = \sum_{k=1}^{d} \prod_{\substack{j=1 \\ j \neq k}}^{d} (z - z_j).
\]
Hence
\[
P'(z_k) = \prod_{\substack{j=1 \\ j \neq k}}^{d} (z_k - z_j) 
= d \prod_{j=1}^{d-1} (z_k - c_j), \quad k=1, \dots, d,
\]
and
\[
P(c_j) = \prod_{k=1}^{d} (c_j - z_k), \quad j=1, \dots, d-1.
\]
This implies
\begin{equation}\label{mie}
\prod_{\substack{k,j=1 \\ j\neq k}}^{d} (z_k - z_j)
= \prod_{k=1}^d P'(z_k) 
= d^d \prod_{k=1}^d \prod_{j=1}^{d-1} (z_k - c_j)
\end{equation}
\[
= (-1)^{d(d-1)} d^d \prod_{j=1}^{d-1} P(c_j).
\]

Let $\C^d$ denote the $d$-dimensional complex vector space and
\[
\Delta^d = \{ \mathbf{z} = (z_1, \dots,z_{k}\cdots, z_d) \in \C^d \;|\; |z_k|<1,\; k=1, \dots, d\}
\]
be the open unit polydisk, with closure
\[
\overline{\Delta^d} = \{ \mathbf{z} = (z_1, \dots,z_{k}\cdots, z_d)\in \C^d \;|\; |z_k| \le 1,\; k=1, \dots, d \}.
\]
Let
\[
\T^d = \{ \mathbf{z}= (z_1, \dots,z_{k}\cdots, z_d) \in \C^d \;|\; |z_k| = 1,\; k=1, \dots, d \}
\]
be the $d$-torus. Denote $\Delta=\Delta^{1}$, $\overline{\Delta}=\overline{\Delta^{1}}$, and ${\mathbb T} ={\mathbb T}^{1}$. Consider the holomorphic map
\[
f(\mathbf{z}) = \prod_{\substack{k,j=1 \\ j\neq k}}^{d} (z_k - z_j) : \C^d \to \C.
\]

\begin{proposition}\label{prop1}
We have
\[
\max_{\mathbf{z} \in \overline{\Delta^d}} \{ |f(\mathbf{z})|\} \le d^d.
\]
The maximum value is attained at $\mathbf{z} = (z_1, \dots, z_d) \in \T^d$ with $z_1, \dots, z_d$ symmetrically distributed on the unit circle.
\end{proposition}

\begin{proof}
Fix $z_2, \dots, z_d$. Then $f(\mathbf{z})$ is holomorphic in $z_1 \in \Delta$. 
By the maximum modulus principle, $|f|$ attains its maximum on $|z_1|=1$. 
Repeating this argument for each $z_k$, we conclude that the maximum is attained on $\T^d$.

Schur's classical inequality \cite{Sh} (see also~\cite[Proposition 2.3]{T}) shows that the maximum value $d^d$ is achieved for symmetrically distributed points on $\T^d$.
\end{proof}

\begin{example}
For the quadratic $P(z) = z^2 + a_1 z + a_0$ with roots $z_1, z_2$, we have
\[
f(\mathbf{z}) = -(z_1 - z_2)^2.
\]
The maximum value of $|f|$ on $\overline{\Delta^2}$ is attained at $z_1=1$, $z_2=-1$, corresponding to $P(z)=z^2-1$. Hence
\[
\max_{\mathbf{z} \in \overline{\Delta^2}} \{|f(\mathbf{z})|\} \le 2^2.
\]
\end{example}

\begin{example}
For $P(z) = z^4 + a_3 z^3 + a_2 z^2 + a_1 z + a_0$ with roots $z_1, \dots, z_4$, we have
\[
f(\mathbf{z}) = \prod_{\substack{k,j=1 \\ k\neq j}}^{4} (z_k - z_j).
\]
The maximum value of $|f|$ on $\overline{\Delta^4}$ is attained at $z_1=1$, $z_2=i$, $z_3=-1$, $z_4=-i$, corresponding to $P(z)=z^4-1$. Hence
\[
\max_{\mathbf{z} \in \overline{\Delta^4}} \{|f(\mathbf{z})|\} \le |1-(-1)|^4 |1-i|^4 |1-(-i)|^4 = 4^4.
\]
\end{example}

Proposition~\ref{prop1} and (\ref{mie}) imply the statement of \cite[Theorem 2]{D}.

\begin{corollary}
If all roots $z_1, \dots, z_d$ lie in $\overline{\Delta}$, then
\[
\left| \prod_{i=1}^{d-1} P(c_i) \right|^{\frac{1}{d-1}} \le 1.
\]
In particular, there exists a critical point $c$ with
\[
|P(c)| \le 1.
\]
\end{corollary}

Now we prove our main result. 

\medskip
\begin{proof}[Proof of Theorem~\ref{main}]
Let $R>1$ be a real number and set
\[
z_k = R e^{2 \pi i \frac{k}{d-2}},\quad k=1, \dots, d-2,\quad 
z_{d-1} = \frac{(-1)^{d-1}}{R^{d-2}}, \quad z_{d}=0. 
\]
We obtain a monic polynomial 
\[
P(z) = z \prod_{k=1}^{d-1} (z - z_k) 
= z^d + a_{d-1} z^{d-1} + \cdots + a_{2} z^{2} + z \in {\mathcal P}_{1}.
\]

Note that
\[
P'(0)= (-1)^{d-1} \prod_{k=1}^{d-1} z_{k} = (-1)^{d-1} d \prod_{k=1}^{d-1} c_{k} = 1.
\]
Thus
\[
\prod_{k=1}^{d-1} |c_k| = \frac{1}{d}.
\]

Hence
\[
|GM_P| 
= \left(\prod_{k=1}^{d-1} \Big| \frac{P(c_k)}{c_k}\Big|\right)^{\frac{1}{d-1}}
= \left( d \prod_{k=1}^{d-1} |P(c_k)| \right)^{\frac{1}{d-1}}.
\]

Using \eqref{mie}, 
\begin{align*}
d^{d} \prod_{i=1}^{d-1} |P(c_{i})| 
&= \prod_{\substack{k,j=1 \\ j\neq k}}^{d} |z_{k}-z_{j}|  \\
&= \Big(\prod_{\substack{k,j=1 \\ j\neq k}}^{d-1} |z_{k}-z_{j}| \Big)\Big(\prod_{k=1}^{d-1} |z_{k}|\Big)^{2} \\
&= \prod_{\substack{k,j=1 \\ j\neq k}}^{d-1} |z_{k}-z_{j}| \\
&= \Big( \prod_{\substack{k,j=1 \\ j\neq k}}^{d-2} |z_{k}-z_{j}|\Big) \Big(\prod_{j=1}^{d-2} |z_{d-1}-z_{j}|\Big)^{2}.
\end{align*}
In the last formula, when $d=3$, we denote
$$
\prod_{\substack{k,j=1 \\ j\neq k}}^{3-2} |z_{k}-z_{j}| =1.
$$

When $d>3$, for any $0\leq k<j\leq d-2$ or $0\leq j<k\leq d-2$,
$$
|z_{k}-z_{j}| \geq R \Big|1-e^{2\pi i\frac{1}{d-2}}\Big|\geq 2R\Big|1-\cos\Big(\frac{2\pi}{d-2}\Big)\Big|.
$$
Thus, for $d=3$
\[
3^{3} \prod_{i=1}^{3-1} |P(c_{i})|
\geq |z_{2}-z_{1}|^{2}
\]
and for $d> 3$, 
\[
d^{d} \prod_{i=1}^{d-1} |P(c_{i})|
\geq  \Big(2 \Big|1-\cos\Big(\frac{2\pi}{d-2}\Big)\Big| R\Big)^{(d-2)(d-3)} \Big(\prod_{k=1}^{d-2} |z_{d-1}-z_{k}|\Big)^{2}.
\]
Since for $1\leq k\leq d-2$,
\[
|z_{d-1}-z_{k}| \ge R-1,
\]
we get that for $d=3$,
\[
3^{3} \prod_{i=1}^{3-1} |P(c_{i})|
\geq (R-1)^{2}
\]
and for $d>3$,
$$
d^{d} \prod_{i=1}^{d-1} |P(c_{i})| \ge \Big(2 \Big|1-\cos\Big(\frac{2\pi}{d-2}\Big)\Big| R\Big)^{(d-2)(d-3)} (R-1)^{2(d-2)}. 
$$
Therefore, for $d=3$,
\[
3 \prod_{i=1}^{3-1} |P(c_{i})|
\geq  \frac{(R-1)^{2}}{3^{2}}
\]
and, for $d>3$,
\[
d \prod_{i=1}^{d-1} |P(c_{i})| \ge \frac{\Big(2 \Big|1-\cos\Big(\frac{2\pi}{d-2}\Big)\Big| R\Big)^{(d-2)(d-3)}}{d^{d-1}} (R-1)^{2(d-2)}.
\]
This implies that for $d=3$,
\[
|GM_{P}| \geq  \frac{R-1}{3}
\]
and for $d>3$, 
\begin{align*}
|GM_{P}| & \ge \left(\frac{\Big(2 \Big|1-\cos\Big(\frac{2\pi}{d-2}\Big)\Big| R\Big)^{(d-2)(d-3)}}{d^{d-1}} (R-1)^{2(d-2)}\right)^{\frac{1}{d-1}} \\
&=
\frac{\Big(2 \Big|1-\cos\Big(\frac{2\pi}{d-2}\Big)\Big| R\Big)^{\frac{(d-2)(d-3)}{d-1}}}{d} (R-1)^{2\frac{d-2}{d-1}}
\end{align*}
Choosing $R$ sufficiently large, we obtain $|GM_P| > L$.  
This completes the proof. 
\end{proof}

When $d=3$, suppose $c_{1}$ and $c_{2}$ are the two critical points of $P$ with $|c_{1}| \leq |c_{2}|$.  Since 
\[
P'(z) = 3 a_{3} (z - c_{1})(z - c_{2})
= 3 a_{3}\big(z^{2} - (c_{1} + c_{2})z + c_{1} c_{2}\big),
\]
we obtain
\[
P(z) = a_{3} \Big(z^{3} - \tfrac{3(c_{1} + c_{2})}{2} z^{2} + 3 c_{1} c_{2} z\Big).
\]

At $z = c_{1}$,
\[
P(c_{1}) = a_{3} \cdot \tfrac{1}{2}(-c_{1} + 3 c_{2}) c_{1}^{2}.
\]
Since $3 a_{3} c_{1} c_{2} = 1$, it follows that
\[
A_{1,P} = \frac{P(c_{1})}{c_{1}}
= \frac{ \tfrac{1}{2}(-c_{1} + 3 c_{2}) c_{1}}{3 c_{1} c_{2}}
= \frac{1}{6}\Big(-\tfrac{c_{1}}{c_{2}} + 3\Big).
\]
Since $|c_{1}/c_{2}| \le 1$,
\[
|A_{1,P}| \le \tfrac{2}{3} = \tfrac{3-1}{3}.
\]
This verifies that Smale's conjecture holds for $d=3$.

However, $|A_{2,P}|$ and $|GM_{P}|$ may be arbitrarily large since at $z = c_{2}$,
\[
P(c_{2}) = a_{3} \cdot \tfrac{1}{2}(-c_{2} + 3 c_{1}) c_{2}^{2}
\]
and since $3 a_{3} c_{1} c_{2} = 1$, we get
\[
A_{2,P} = \frac{P(c_{2})}{c_{2}}
= \frac{ \tfrac{1}{2}(-c_{2} + 3 c_{1}) c_{2}}{3 c_{1} c_{2}}
= \frac{1}{6}\Big(-\tfrac{c_{2}}{c_{1}} + 3\Big).
\]
Thus
\[
|A_{2,P}| \;\ge\; \tfrac{1}{6}\Big(\Big|\frac{c_{2}}{c_{1}}\Big| - 3\Big),
\]
which may be arbitrarily large.
The geometric mean value is
\[
GM_{P}
= \frac{1}{6} \sqrt{\Big(-\tfrac{c_{1}}{c_{2}} + 3\Big)\Big(-\tfrac{c_{2}}{c_{1}} + 3\Big)}
= \frac{1}{6} \sqrt{8 - 3\Big(\tfrac{c_{1}}{c_{2}} + \tfrac{c_{2}}{c_{1}}\Big)}.
\]
Since
\[
\Big|8 - \tfrac{3 c_{1}}{c_{2}}\Big|
\le 8 + \Big|\tfrac{3 c_{1}}{c_{2}}\Big|
\le 11,
\]
we have
\[
|GM_{P}| \;\ge\; \frac{1}{6} \sqrt{ \Big|\tfrac{3 c_{2}}{c_{1}}\Big| - \Big|8 - \tfrac{3 c_{1}}{c_{2}}\Big| }
\;\ge\; \frac{1}{6} \sqrt{3\Big|\frac{c_{2}}{c_{1}}\Big| - 11},
\]
which may again be arbitrarily large.

\section{\bf \large On Smale's Conjecture}~\label{sec:smvc}

Fixing $d>2$, let  
\[
P(z) = a_{d} z^{d} + a_{d-1} z^{d-1} + \cdots + a_{2} z^{2} + z \in \mathcal{P}_{1}, \;\;a_{d}\not=0.
\]
Suppose $z_{1},\dots,z_{d-1}$ are all the nonzero roots of $P(z)=0$ (numbered with multiplicities taken into account) and suppose $c_{1},\dots,c_{d-1}$ are all the critical points of $P$ (numbered with multiplicities taken into account). Then
\[
P(z) = a_{d} z \prod_{k=1}^{d-1} (z - z_{k}),
\quad 
P'(z) = a_{d} d \prod_{k=1}^{d-1} (z - c_{k}).
\]
Since $P'(0)=1$,
$$
(-1)^{d-1} a_{d}\prod_{k=1}^{d-1} z_{k} =(-1)^{d-1} a_{d}d \prod_{k=1}^{d-1} c_{k} =1.
$$
This implies that 
\begin{equation}~\label{rc}
\prod_{k=1}^{d-1} z_{k} =\frac{(-1)^{d-1}}{a_{d}}\quad \hbox{and}\quad
\prod_{k=1}^{d-1} c_{k} =\frac{(-1)^{d-1}}{a_{d} d}.
\end{equation}
Consider the monic polynomial
\[
Q(z) = \frac{P(z)}{a_{d}} 
= z^{d} + \frac{a_{d-1}}{a_{d}} z^{d-1} + \cdots + \frac{a_{2}}{a_{d}} z^{2} + \frac{1}{a_{d}} z.
\]
Then 
$$
Q'(0) =\frac{P'(0)}{a_{d}} =\frac{1}{a_{d}}=(-1)^{d-1} \prod_{k=1}^{d-1} z_{k} =(-1)^{d-1} d \prod_{k=1}^{d-1} c_{k}.
$$
Since 
\[
A_{k,P} = \frac{P(c_{k})}{c_{k}}
= \frac{Q(c_{k})}{c_{k}} \cdot \frac{1}{Q'(0)}=a_{d} \frac{Q(c_{k})}{c_{k}}, 
\quad k = 1,\dots,d-1,
\]
we have that 
\[
GM_{P} = \Bigg( \prod_{k=1}^{d-1} A_{k,P} \Bigg)^{\!\!\frac{1}{d-1}}
= \Big( a_{d}^{d-1} \prod_{k=1}^{d-1} \frac{Q(c_{k})}{c_{k}} \Big)^{\frac{1}{d-1}}
= a_{d} \Bigg( \prod_{k=1}^{d-1} \frac{Q(c_{k})}{c_{k}} \Bigg)^{\frac{1}{d-1}}.
\]
Thus,
$$
|GM_{P}| =|a_{d}| \Big(|a_{d}|d \prod_{k=1}^{d-1} |Q(c_{k})| \Big)^{\frac{1}{d-1}}= |a_{d}|^{\frac{d}{d-1}} d^{\frac{1}{d-1}} \Big( \prod_{k=1}^{d-1} |Q(c_{k})|\Big)^{\frac{1}{d-1}}.
$$

\medskip
\begin{theorem}\label{unit}
If $|z_{k}| \leq 1$ for all $k=1,\dots,d-1$, then
\[
|GM_{P}| \leq |a_{d}|^{\frac{d-2}{d-1}} \cdot \frac{d-1}{d} \;\leq\; |a_{d}| \cdot \frac{d-1}{d}.
\]
\end{theorem}

\begin{proof}
Equality~(\ref{mie}) gives
\begin{align*}
d^{d} \prod_{i=1}^{d-1} |Q(c_{i})|
&= \prod_{\substack{k,j=1 \\ k \neq j}}^{d} |z_{k} - z_{j}|\\
&= \Bigg( \prod_{\substack{k,j=1 \\ k \neq j}}^{d-1} |z_{k} - z_{j}| \Bigg) 
\left( \prod_{k=1}^{d-1} |z_{k}| \right)^{2}\\
&= |Q'(0)|^{2} \prod_{\substack{k,j=1 \\ k \neq j}}^{d-1} |z_{k} - z_{j}|\\
&= \frac{1}{|a_{d}|^{2}} \prod_{\substack{k,j=1 \\ k \neq j}}^{d-1} |z_{k} - z_{j}|.
\end{align*}

By Proposition~\ref{prop1},
\[
\prod_{i=1}^{d-1} |Q(c_{i})| 
= \frac{1}{|a_{d}|^{2} d^{d}} \prod_{k \neq j}^{d-1} |z_{k} - z_{j}| 
\;\leq\; \frac{(d-1)^{d-1}}{|a_{d}|^{2}d^{d}}.
\]
Hence
\begin{align*}
|GM_{P}| &=|a_{d}|^{\frac{d}{d-1}} d^{\frac{1}{d-1}} \Big(\prod_{i=1}^{d-1} |Q(c_{i})|\Big)^{\frac{1}{d-1}} \leq |a_{d}|^{\frac{d}{d-1}} d^{\frac{1}{d-1}} \Big( \frac{(d-1)^{d-1}}{|a_{d}|^{2}d^{d}}\Big)^{\frac{1}{d-1}}\\
&\leq\; \frac{d-1}{d} \cdot |a_{d}|^{\frac{d-2}{d-1}}.
\end{align*}
Since
\[
|a_{d}| = \frac{1}{\prod_{k=1}^{d-1} |z_{k}|} \geq 1,
\]
we have $|a_{d}|^{(d-2)/(d-1)} \leq |a_{d}|$.  
Thus,
\[
|GM_{P}| \leq |a_{d}| \cdot \frac{d-1}{d}.
\]
\end{proof}

When $P\in {\mathcal P}_{1}$ is a monic polynomial, that is, $a_{d}=1$, and when $|z_{k}|\leq 1$, $k=1$, $\cdots$, $d-1$, then $|z_{k}|=1$, $k=1$, $\cdots$, $d-1$, since $P'(0)=|\prod_{k=1}^{d-1}z_{k}|=1$.  
Thus, the result first established in~\cite[Proposition~2.3]{T} becomes a corollary of Theorem~\ref{unit}.

\medskip
\begin{corollary}~\label{eqroot}
Smale's conjecture holds for all $P \in \mathcal{P}_{1}$ 
such that all the nonzero roots of $P(z) = 0$ have the same norm. 
\end{corollary}

\begin{proof}
Suppose 
$$
|z_{1}|=\cdots =|z_{d-1}| =\lambda.
$$
Since 
\[
P(z) = a_{d} z \prod_{k=1}^{d-1} (z - z_{k})
\] 
and since $P'(0)=1$, we have that $|a_{d}|\lambda^{d-1}=1$.
This implies that 
$$
\lambda =|a_{d}|^{-\frac{1}{d-1}}
$$
and consider 
\[
P_{\lambda}(z) = \frac{1}{\lambda} P(\lambda z) 
= a_{d} \lambda^{d-1} z \prod_{k=1}^{d-1} \Big(z - \frac{z_{k}}{\lambda}\Big)\in {\mathcal P}_{1}.
\]
Then $P_{\lambda}\in {\mathcal P}_{1}$ is a monic polynomial. The roots of $P_{\lambda}(z)=0$ are $0$ and 
$$
z_{i, \lambda}=\frac{z_{i}}{\lambda}, \quad 1\leq i\leq d-1.
$$
Theorem~\ref{unit} implies 
\[
|GM_{P_{\lambda}}| \leq \frac{d-1}{d}.
\]

Since 
$$
A_{k, P} =A_{k, P_{\lambda}}, \quad k=1,\cdots, d-1,
$$
\[
|GM_{P}| =|GM_{P_{\lambda}}| \leq \frac{d-1}{d}.
\]
This implies that 
$$
\min\{|A_{k, P}|\; ;\; 1\leq k\leq d-1 \} \leq \frac{d-1}{d}.
$$
This completes the proof.
\end{proof}

\medskip
The following inequality is adapted from~\cite{D}.

\begin{lemma}\label{eq0}
For any 
\[
P(z) = a_{d} z^{d} + \cdots + z \in \mathcal{P}_{1},
\]
we have
\[
\min_{1 \le k \le d-1} \{|P(c_{k})|\}
\;\leq\; 
\frac{d-1}{d} \, |d a_{d}|^{-\tfrac{1}{d-1}}.
\]
\end{lemma}

Using this lemma, we can confirm Smale's conjecture under suitable assumptions on the critical points. 
For $d>2$ and 
\[
P(z) = a_{d} z^{d} + \cdots + a_{2} z^{2} + z \in \mathcal{P}_{1},
\]
let $c_{1}$, $\cdots$, $c_{d-1}$ be its critical points (numbered with multiplicities taken into account) and $k$ be an index such that
$$
|P(c_{k})|= \min_{1\leq l\leq d-1}\{ |P(c_{l})|\}.
$$

\medskip
\begin{theorem}~\label{eqcritsymm}
Smale's conjecture holds for $P \in {\mathcal P}_{1}$ with the absolute value of the geometric mean of critical points less than or equal to $|c_{k}|$, that is, 
\[
\Big( \prod_{l=1}^{d-1} |c_{l}|\Big)^{\frac{1}{d-1}} \leq |c_{k}|.
\]
\end{theorem} 

\begin{proof}
Since $P'(z) =a_{d} d \prod_{i=1}^{d-1} (z-c_{i})$,
\[
1=|P'(0)| = |a_{d}d| \prod_{k=1}^{d-1} |c_{k}| \leq |a_{d}d| |c_{k}|^{d-1}.
\]
This implies that $|a_{d}d|^{-1} \leq  |c_{k}|^{(d-1)}$.  
Applying Lemma~\ref{eq0}, we have that
\[
|P(c_{k})| \;\le\; \frac{d-1}{d} |a_{d} d|^{-\frac{1}{d-1}}\leq \frac{d-1}{d}|c_{k}|. 
\]
Hence 
\[
|A_{k,P}| =\frac{|P(c_{k})|}{|c_{k}|} \leq \frac{d-1}{d}.
\]
This implies that 
$$
\min\{ |A_{i, P}|\; ;\; 1\leq i\leq d-1\} \leq \frac{d-1}{d}.
$$
\end{proof}

A direct consequence of Theorem~\ref{eqcritsymm} is 
that Smale's conjecture holds when all critical points have the same norm. 

\medskip
\begin{corollary}\label{eq3}
Smale's conjecture holds for all $P \in \mathcal{P}_{1}$ whose critical points have the same norm.
\end{corollary}

Suppose $P\in {\mathcal P}_{1}$. Since (\ref{rc}) and since 
\[
P'(z)=d a_d\prod_{j=1}^{d-1}(z-c_j),
\]
we have that 
\[
P'(z_k)=d a_d\prod_{j=1}^{d-1}(z_k-c_j) =\frac{(-1)^{d-1}}{\prod_{j=1}^{d-1}c_j}
  \prod_{j=1}^{d-1}(z_k-c_j) =\prod_{j=1}^{d-1}\left(1-\frac{z_k}{c_j}\right).
\]
Also, for 
\[
P(z)=a_d z\prod_{k=1}^{d-1}(z-z_k),
\]
\[
\begin{aligned}
\frac{P(c_j)}{c_j} =a_d\prod_{k=1}^{d-1}(c_j-z_k)
=\frac{(-1)^{d-1}}{\prod_{k=1}^{d-1}z_k}
  \prod_{k=1}^{d-1}(c_j-z_k)=\prod_{k=1}^{d-1}\left(1-\frac{c_j}{z_k}\right).
\end{aligned}
\]
Moreover
$$
GM_{P} =\Big(\prod_{k=1}^{d-1} \frac{P(c_j)}{c_j}\Big)^{\frac{1}{d-1}}  =\frac{1}{d} \Big( \prod_{k=1}^{d-1}P'(z_{k})\Big)^{\frac{1}{d-1}}
$$

For each \(1\leq k\leq d-1\), define the Lagrange interpolation identity
\[
L_k(z)=\prod_{\substack{j=1\\j\ne k}}^{d-1}
       \frac{z-z_j}{z_k-z_j}, 
\]
\[
L_k(z_k)=1,\;\;
L_k(z_j)=0, \;\; 1\le j\ne k\le d-1, \;\; \hbox{and} \;\;
L_k(0)=\prod_{\substack{j=1\\j\ne k}}^{d-1}
       \frac{-z_j}{z_k-z_j}.
\]
Define
\[
F(z)=\sum_{k=1}^{d-1}L_k(z).
\]
Since \(\deg F\le d-2\) and \(F(z_k)=1\) for
\(k=1,\ldots,d-1\), it follows that
\[
F(z)\equiv 1.
\]
In particular,
\[
\sum_{k=1}^{d-1}L_k(0)=1.
\]
Since 
\[
P'(z_k)=a_d z_k
\prod_{\substack{j=1\\j\ne k}}^{d-1}(z_k-z_j)
=\prod_{\substack{j=1\\j\ne k}}^{d-1}
  \left(1-\frac{z_k}{z_j}\right),
\]
we have that 
\[
\frac{1}{P'(z_k)}
=\prod_{\substack{j=1\\j\ne k}}^{d-1}
  \frac{z_j}{z_j-z_k}
  =(-1)^{d-2}L_k(0).
\]
Thus implies that 
\[
\sum_{k=1}^{d-1}\frac{1}{P'(z_k)}
=(-1)^{d-2}\sum_{k=1}^{d-1}L_k(0)
=(-1)^{d-2},
\]
and 
$$
\sum_{k=1}^{d-1}\frac{1}{|P'(z_k)|} \geq 1.
$$
Thus,
$$
\max \Big\{ \frac{1}{|P'(z_{k})|}\; ;\; 1\leq k\leq d-1 \Big\} \geq \frac{1}{d-1}.
$$
Therefore, we proved the following theorem in the spirit of Smale's conjecture. 

\medskip
\begin{theorem}
For any $P\in {\mathcal P}_{1}$, 
\[
\min \{ |P'(z_k)|\; ;\; 1\le k\le d-1\} \le d-1,
\]
or equivalently,
\[
\min \Big\{ \left|\frac{P'(z_k)}{d}\right|\; ;\; 1\le k\le d-1\Big\} \le \frac{d-1}{d}.
\]
\end{theorem}

Let \(P(z)=zR(z)\). Then 
\[
P'(z)=R(z)+zR'(z), \quad R(0)=1,
\]
and
\[
R(c_i)+c_iR'(c_i)=0,\;\; \frac{P(c_i)}{c_i}=R(c_i),\quad i=1,\ldots,d-1.
\]
Therefore, we have that 
$$
GM_{P} =\Big( \prod_{i=1}^{d-1} R(c_{i})\Big)^{\frac{1}{d-1}}
$$
Since 
\[
P'(z)=\prod_{j=1}^{d-1}\left(1-\frac{z}{c_j}\right),
\]
\[
P(z)=\int_0^z\prod_{j=1}^{d-1}
\left(1-\frac{t}{c_j}\right)\,dt
\]
and
\[
R(c_i)=\frac{1}{c_i}\int_0^{c_i}
\prod_{j=1}^{d-1}\left(1-\frac{t}{c_j}\right)\,dt.
\]
The Lagrange interpolation gives
\[
R(z)=\sum_{i=1}^{d-1}R(c_i)
\prod_{\substack{j=1\\j\ne i}}^{d-1}
\frac{z-c_j}{c_i-c_j}.
\]
This implies that 
\[
 R(0)=\sum_{i=1}^{d-1}
R(c_i)\prod_{\substack{j=1\\j\ne i}}^{d-1}
\frac{-c_j}{c_i-c_j}=1.
\]
Thus,
\[
\sum_{i=1}^{d-1}
\Big| R(c_i)\prod_{\substack{j=1\\j\ne i}}^{d-1}
\frac{-c_j}{c_i-c_j}\Big|\geq 1
\]
and
$$
\max \Big\{ \Big| R(c_i)\prod_{\substack{j=1\\j\ne i}}^{d-1}
\frac{-c_j}{c_i-c_j}\Big|\; ;\; 1\leq i\leq d-1\Big\} \geq \frac{1}{d-1}.
$$
Thus, we have another theorem in the spirit of Smale's conjecture. 

\medskip
\begin{theorem}
For $P\in {\mathcal P}_{1}$, 
\[
\min 
\Big\{\frac{\prod_{j=1, j\ne i}^{d-1}
 \Big|1-\frac{c_i}{c_j}\Big|}{|A_{i,P}|}\; ;\; 1\leq i\leq d-1\Big\} \leq d-1
\]
or equivalently,
\[
\min 
\Big\{\frac{\prod_{j=1, j\ne i}^{d-1}
 \Big|1-\frac{c_i}{c_j}\Big|}{d|A_{i,P}|}\; ;\; 1\leq i\leq d-1\Big\} \leq \frac{d-1}{d}
\]
\end{theorem}





\section{\bf \large Smale's Mean Values and Fekete Points}~\label{sec:smvfp}

From (\ref{norm}), to study Problem~\ref{scp}, we only need to study the subspace $\mathcal{P}_{2}={\mathcal P}_{2}(d)$ consisting of all
$P\in {\mathcal P}$ such that
$$
P'(z) =(z-c_{1})\cdots (z-c_{d-1}), \;\;
0<|c_{1}|\leq |c_{2}|\leq \cdots \leq |c_{d-1}|=1. 
$$
For $P\in {\mathcal P}_{2}$, 
$$
P'(0)=(-1)^{d-1}c_{1}\cdots c_{d-1}\not=0
$$
and
$$
P'(z)=z^{d-1}+ \sum_{k=1}^{d-1} (-1)^{k} \Big( \sum_{1\leq j_{1}<j_{2}<\cdots<j_{k}\leq d-1} c_{j_{1}}\cdots c_{j_{k}}\Big) z^{d-1-k}.
$$
Then 
$$
P(z)=\frac{1}{d} z^{d} +\sum_{k=1}^{d-1} \frac{(-1)^{k}}{d-k} \Big( \sum_{1\leq j_{1}<j_{2}<\cdots<j_{k}\leq d-1} c_{j_{1}}\cdots c_{j_{k}}\Big) z^{d-k}.
$$
Suppose  
\begin{align*}
P(z) & =\frac{1}{d} z (z-z_{1})\cdots (z-z_{d-1})\\
    & =\frac{1}{d} z\Big(z^{d-1} +\sum_{k=1}^{d-1} (-1)^{k} \Big( \sum_{1\leq j_{1}<j_{2}<\cdots<j_{k}\leq d-1} z_{j_{1}}\cdots z_{j_{k}} \Big) z^{d-1-k}\Big)\\
    &=\frac{1}{d} z^{d} + \sum_{k=1}^{d-1} \frac{(-1)^{k}}{d}\Big( \sum_{1\leq j_{1}<j_{2}<\cdots<j_{k}\leq d-1} z_{j_{1}}\cdots z_{j_{k}} \Big) z^{d-k}.
\end{align*}
Comparing with the previous formula, we get, for $1\leq k\leq d-1$,
$$
\frac{1}{d} \Big(\sum_{1\leq j_{1}<j_{2}<\cdots<j_{k}\leq d-1} z_{j_{1}}\cdots z_{j_{k}}\Big) =
\frac{1}{d-k} \Big( \sum_{1\leq j_{1}<j_{2}<\cdots<j_{k}\leq d-1} c_{j_{1}}\cdots c_{j_{k}}\Big).
$$
Thus, 
$$
P(z)=\frac{1}{d} z^{d} + \sum_{k=1}^{d-1} \frac{(-1)^{k}}{d-k}\Big( \sum_{j_{1}<j_{2}<\cdots<j_{k}}^{d-1} c_{j_{1}}\cdots c_{j_{k}}\Big) z^{d-k}.
$$
This implies that 
\begin{align*}
A_{l,P} & = \frac{P(c_{l})}{c_{l}}\frac{1}{P'(0)}\\
        & =\Big( \frac{1}{d} c_{l}^{d} + \sum_{k=1}^{d-1}\Big(\frac{(-1)^{k}}{d-k}\sum_{1\leq j_{1}<j_{2}<\cdots<j_{k}\leq d-1} c_{j_{1}}\cdots c_{j_{k}}\Big) c_{l}^{d-k} \Big)\frac{(-1)^{d-1}}{c_{1}\cdots c_{l}^{2}\cdots c_{d-1}}\\
& =\frac{1}{2} +\sum_{k=0}^{d-2} \frac{(-1)^{k}}{(k+2)(k+1)}
\sum_{1\leq j_{1}<j_{2}<\cdots<j_{k+1}\leq d-1, j_{s}\not=l} \frac{c_{l}}{c_{j_{1}}} \cdots \frac{c_{l}}{c_{j_{k+1}}}.
\end{align*}
Rewriting this, we have for $1\leq l\leq d-1$, 
\begin{equation}~\label{feq}
A_{l,P}=\frac{1}{2} +\sum_{k=0}^{d-2} \frac{(-1)^{k}}{(k+2)(k+1)}
\sum_{1\leq j_{1}<j_{2}<\cdots<j_{k+1}\leq d-1, j_{s}\not=l} \frac{c_{l}}{c_{j_{1}}} \cdots \frac{c_{l}}{c_{j_{k+1}}}
\end{equation}

Consider 
$$
(\overline{\Delta^{d-1}})^{*} =\{ {\bf c} =(c_{1}, \cdots, c_{d-1})\in \C^{d-1}\;|\; 0<|c_{1}|\leq \cdots \leq |c_{d-1}|\leq 1\}.
$$
Equation (\ref{feq}) defines a several-complex-variable function on $(\overline{\Delta^{d-1}})^{*}$,  
$$
F(c_{1}, \cdots, c_{l}, \cdots, c_{d-1})=(A_{1,P}, \cdots, A_{l,P},\cdots, A_{d-1,P}).
$$
Now Problem~\ref{scp} is equivalent to the following problem. 

\medskip
\begin{problem}~\label{scp1}
Find a smallest number $M>0$ such that the image $F((\overline{\Delta^{d-1}})^{*})$ is contained in the closed domain
$$
D^{d-1} =\{ {\bf w} =(w_{1}, \cdots, w_{l}, \cdots, w_{d-1})\;|\; \hbox{at least one $l$ with $|w_{l}|\leq M$}\}
$$
in ${\mathbb C}^{d-1}$. 
\end{problem}

The advantage of the existence of the several-complex-variable function $F$ is that we can use the computer to test Problem~\ref{scp}. 
Furthermore, we can formulate a more general problem.
Consider 
\begin{align*}
&f_{l, d, m} (c_{1}, \cdots, c_{d-1})=\frac{(-1)^{m-d}}{(m-d+1)(m-d+2)}\\
&+\sum_{k=0}^{d-2} \frac{(-1)^{m-d+k}}{(m-d+k+2)(m-d+k+3)} \sum_{1\leq j_{1}<j_{2}<\cdots<j_{d-k}\leq d-1, j_{s}\not=l} \frac{c_{l}}{c_{j_{1}}} \cdots \frac{c_{l}}{c_{j_{d-k}}}.
\end{align*}
Let 
$$
L_{m,d} =\sum_{k=1}^{d-1}\frac{1}{(m-k)(m-k+1)}= \frac{1}{m-d+1}-\frac{1}{m}.
$$
Define a several-complex-variable function on $(\overline{\Delta^{d-1}})^{*}$,
$$
F_{m,d} (c_{1}, \cdots, c_{d-1})= (f_{1,m,d}, \cdots, f_{d-1,m,d}).
$$

\medskip
\begin{problem}~\label{gsc}
Can we find a number $M>0$ such that the image $F_{m,d}((\overline{\Delta^{d-1}})^{*})$ is contained in 
$$
D^{d-1}_{m}=\{ {\bf w} =(w_{1}, \cdots, w_{d-1})\;|\; \hbox{at least one $1\leq l\leq d-1$ with $|w_{l}|\leq M$}\}
$$
in ${\mathbb C}^{d-1}$.
\end{problem}

The conjecture is that $M\leq 4 L_{m,d}$. Problems~\ref{scp1} and~\ref{gsc} connect Problem~\ref{scp} to Fekete points, a specific set of points on a compact set in the complex plane that maximizes a particular product of distances, providing a near-optimal distribution for interpolation and approximation.  

Given a point ${\bf c}=(c_{1}, \cdots, c_{d-1})\in \C^{d-1}$, consider the function
\begin{align*}
f_{i} (z) &=\prod_{1\leq j\leq d-1, j\not=i} \Big( 1-\frac{c_{i}}{c_{j}} z\Big)\\
&= 1+ \sum_{k=1}^{d-2} \sum_{1\leq j_{1}<j_{2}<\cdots<j_{k}\leq d-1, j_{s}\not=i} (-1)^{k} \frac{c_{i}}{c_{j_{1}}} \cdots \frac{c_{i}}{c_{j_{k}}} z^{k}
\end{align*}
In the theory of Fekete points, one studies the supremum and the infimum of 
$f_{1}(1), \cdots, f_{d-1}(1)$ where ${\bf c}$ runs over $(\overline{\Delta^{d-1}})^{*}$. Furthermore,
consider
\begin{align*}
f(1) &=\prod_{i=1}^{d-1} f_{i}(1)= \prod_{i=1}^{d-1} \prod_{1\leq j\leq d-1, j\not=i} \Big(1-\frac{c_{i}}{c_{j}}\Big)\\
&=\frac{1}{(c_{1}\cdots c_{d-1})^{d-2}} \prod_{1\leq j\not= i\leq d-1} (c_{j}-c_{i})\\
&=\frac{(-1)^{d-2}}{(c_{1}\cdots c_{d-1})^{d-2}} \prod_{1\leq j<l\leq d-1} (c_{j}-c_{i})^{2}.
\end{align*}
One studies the supremum and the infimum of $f(1)$ where ${\bf c}$ runs over $(\overline{\Delta^{d-1}})^{*}$. 

We now relate the Fekete products $\{ f_{i} (z)\}_{i=1}^{d-1}$ with the Smale's mean values $\{ A_{i, P}\}_{i=1}^{d-1}$ for $P\in {\mathcal P}_{2}$ with critical points ${\bf c} =(c_{1}, \cdots, c_{d-1})$ and 
$
0<|c_{1}| \leq \cdots \leq |c_{d-1}|=1.
$

Let
\begin{align*}
F_{i} (z) &=\int_{0}^{z} f_{i}(\xi) d\xi\\
&=\int_{0}^{z} \Big(1+ \sum_{k=1}^{d-2} \sum_{1\leq j_{1}<j_{2}<\cdots<j_{k}\leq d-1, j_{s}\not=i, s=1, \cdots, k} (-1)^{k} \frac{c_{i}}{c_{j_{1}}} \cdots \frac{c_{i}}{c_{j_{k}}} \xi^{k} \; d\xi\Big)\\
&= z+ \sum_{k=1}^{d-2} \frac{(-1)^{k}}{k+1} \sum_{1\leq j_{1}<j_{2}<\cdots<j_{k}\leq d-1, j_{s}\not=i, s=1, \cdots, k}  \frac{c_{i}}{c_{j_{1}}} \cdots \frac{c_{i}}{c_{j_{k}}} z^{k+1}
\end{align*}
and 
\begin{align*}
G_{i} (z) &=\int_{0}^{z} F_{i}(\xi) d\xi\\
&=\int_{0}^{z} \Big(z+ \sum_{k=1}^{d-2} \frac{(-1)^{k}}{k+1} \sum_{1\leq j_{1}<j_{2}<\cdots<j_{k}\leq d-1, j_{s}\not=i, s=1, \cdots, k} \frac{c_{i}}{c_{j_{1}}} \cdots \frac{c_{i}}{c_{j_{k}}} \xi^{k+1} \; d\xi\Big)\\
&= \frac{1}{2} z^{2}+ \sum_{k=1}^{d-2} \frac{(-1)^{k}}{(k+1)(k+2)} \sum_{1\leq j_{1}<j_{2}<\cdots<j_{k}\leq d-1, j_{s}\not=i, s=1, \cdots, k}  \frac{c_{i}}{c_{j_{1}}} \cdots \frac{c_{i}}{c_{j_{k}}} z^{k+2}
\end{align*}
We can see from (\ref{feq}), 
$$
A_{i, P}=G_{i} (1),  \quad 1\leq i\leq d-1.
$$
This makes the study of Problem~\ref{scp} into the study of Fekete points, that is,

\medskip
\begin{problem}~\label{scf}
Find the minimal number $M>0$ such that 
$$
\min_{1\leq i\leq d-1}\{ |G_{i}(1)|\} \leq M.
$$
\end{problem}

Note that 
$$
\frac{d-1}{d}= 1-\frac{1}{d} = \frac{1}{2} +\sum_{k=1}^{d-2} \Big| \frac{(-1)^{k}}{(k+1)(k+2)}\Big|.
$$
Thus, to study Fekete points, we prove the following theorem, which is in the spirit of Smale's conjecture.

\medskip
\begin{theorem}
$$
\min \{ |f_{i} (1)|\; ;\; 1\leq i \leq d-1\} \leq d-1=1+\sum_{k=1}^{d-2} |(-1)^{k}|.
$$
\end{theorem}

\begin{proof}
If there are two $i\not= j$ such that $c_{i}=c_{j}$, then $f_{i} (1)=0$. Thus, the inequality holds.

Assume all $c_{1}$, $\cdots$, $c_{d-1}$ are distinct numbers. 
Consider the $i$-th Lagrange basis polynomial,
$$
L_i(z) = \prod_{1\leq j\leq d-1, j \neq i} \frac{z - c_j}{c_i - c_j}.
$$
It is a degree $d-2$ polynomial such that 
$$
L_{i} (c_{i}) =1, \quad L_{i}(c_{j}) =0, \;\; 1\leq j\neq i \leq d-1.  
$$
This implies that 
$$
F(z)= \sum_{i=1}^{d-1} L_{i} (z) \equiv 1
$$
since $F(z)$ is a degree $d-2$ polynomial and $F(c_{i})=1$ for $1\leq i\leq d-1$. 

Since 
$$
\frac{1}{f_{i}(1)} =\frac{1}{\prod_{1\leq j\leq d-1, j\not=i} \Big( 1-\frac{c_{i}}{c_{j}} \Big)}=\prod_{1\leq j\leq d-1, j\not=i}\frac{c_{j}}{c_{j}-c_{i}}=(-1)^{d-2} L_{i}(0),
$$
$$
\sum_{i=1}^{d-1} \frac{1}{f_{i}(1)} =(-1)^{d-2} \sum_{i=1}^{d-1} L_{i}(0) =(-1)^{d-2}.
$$

Now we prove the theorem by contradiction. Suppose 
$$
\min\{ |f_{i}(1)|\; ;\; 1\leq i\leq d-1 \} > d-1.
$$
Then 
$$
\frac{1}{|f_{i}(1)|} < \frac{1}{d-1}, \quad \forall\; 1\leq i\leq d-1.
$$
Thus, we have 
$$
1=\Big| \sum_{i=1}^{d-1} \frac{1}{f_{i}(1)} \Big|\leq 
\sum_{i=1}^{d-1} \frac{1}{|f_{i}(1)|} < (d-1)\frac{1}{d-1}=1,
$$
which is impossible. The contradiction proves the theorem.
\end{proof}

We would like to study that 

\medskip
\begin{conjecture}
$$
\min_{1\leq i\leq d-1}\{ |F_{i}(1)|\} \leq 1 +\sum_{k=1}^{d-2} \Big| \frac{(-1)^{k}}{k+1}\Big|=\sum_{k=0}^{d-2} \frac{1}{k+1}.
$$
\end{conjecture}



\bibliographystyle{amsalpha}




 \end{document}